\documentclass[11pt]{amsart}
\usepackage{fullpage, amsmath , amsfonts, amscd, algorithmic, algorithm}

\theoremstyle{plain}
  \newtheorem{theorem}{Theorem}
  \newtheorem{proposition}[theorem]{Proposition}
  \newtheorem{lemma}[theorem]{Lemma}

\theoremstyle{definition}
  \newtheorem{definition}[theorem]{Definition}
  
\theoremstyle{remark}

\newcommand{\field}[1]{\mathbb{#1}}

\def\I{\mathcal{I}}

\def\M{\mathcal{M}}

\def\Z{\field{Z}}

\newcommand{\newword}[1]{\textbf{\emph{#1}}}

\title{Contraction and restriction of positroids in terms of decorated permutations}
\author{Suho Oh}

\begin{document}

\begin{abstract}
A positroid is a matroid defined by Postnikov to study the cells in the nonnegative part of the Grassmannian. They are in bijection with decorated permutations. We show a way to explain contraction and restriction of positroids in terms of decorated permutations. 
\end{abstract}

\maketitle

\section{Introduction}

Matroids are useful objects that generalize certain aspects appearing in hyperplane arrangements, graphs, matrices, etc. Contraction, restriction and dual are operations used often to study matroids. The main focus of this paper, is to study such matroid operations on a special class of matroids, called \newword{positroids}.


Positroids were defined by Postnikov in \cite{postnikov-2006}, to study the nonnegative part of the Grassmannian. This combinatorial object turned out to have remarkable connections to other fields  \cite{KL-KP}, \cite{scatter}. One of the interesting combinatorial properties of positroids is that they can be indexed by combinatorial objects called \newword{Grassmann necklaces} and \newword{decorated permutations}. In this paper, we will study how the matroid operations act on Grassmann necklaces and decorated permutations.

\section{Matroid operations}

In this section, we will review the basics of matroids and the matroid operations. We will be viewing matroids using \newword{bases}. For a more detailed review, we would like to refer the readers to \cite{oxley}.

Let $[n]$ stand for the set $\{1,\ldots,n\}$ and let ${[n] \choose k}$ stand for the collection of all $k$-element subsets of $[n]$. A \newword{matroid} $\M$ is a sub-collection of ${[n] \choose k}$ that satisfies \newword{basis exchange axiom}: for any $I,J \in \M$, pick an arbitrary element $i \in I \setminus J$. Then there exists $j \in J \setminus I$ such that $I \setminus \{i\} \cup \{j\}$ is in $\M$. Matroids appear in numerous different structures, such as hyperplane arrangements, graphs and matrices. In particular, given a $k$-by-$n$ real matrix, by collecting all $k$-element subsets of columns where the minor is nonzero, one can get a matroid.

Let $\M$ be a matroid. The \newword{contraction} of $\M$ by $T \subseteq [n]$ is defined as
$$ \M/T = \{I \setminus T | T \subset I \in \M \} .$$
 The \newword{restriction} of $\M$ to $T \subset [n]$ is defined as
$$ \M|_T = \{I \cap T | |I \cap T| \text{ is maximal among all } I \in \M \}. $$

 The \newword{dual} of $\M$ is defined as
 $$ \M^* = \{I^c | I \in \M\}. $$

Contraction, restriction or dual of a matroid is also a matroid. It is a well known fact that the dual operation of contraction is restriction. 

\begin{theorem}[\cite{oxley}]
If $\M$ is a matroid, $\M^*|_{[n] \setminus T} = (\M/T)^*$.
\end{theorem}

In other words, contraction can be emulated by using restriction and dual.



\section{Positroids, Grassmann necklaces and decorated permutations}

In this section, we will review the definition of positroids and their connection to Grassmann necklaces and decorated permutations. For a more detailed introduction, we would like to refer the readers to \cite{postnikov-2006}.

Positroid is a matroid that can be realized using a real matrix with nonnegative maximal minors. In other words, $\M \subseteq {[n] \choose k}$ is a positroid if and only if there exists a $k$-by-$n$ real matrix such that the maximal minors coming from column set $I \in \M$ are strictly positive, and all other maximal minors are zero.

One of the interesting properties of positroids is that each positroid can be indexed with a \newword{decorated permutation}.

\begin{definition}[\cite{postnikov-2006}, Definition 13.3]
 A decorated permutation $\pi^{:} = (\pi, col)$ is a permutation $\pi \in S_n$ together with a coloring function $col$ from the set of fixed points $\{i | \pi(i) = i\}$ to $\{1,-1\}$. That is, a decorated permutation is a permutation with fixed points colored in two colors.
\end{definition}

To see the connection between decorated permutations and positroids, we use an intermediary combinatorial object called \newword{Grassmann necklace}.

\begin{definition}[\cite{postnikov-2006}, Definition 16.1]
A \textit{Grassmann necklace} is a sequence $\I = (I_1, \cdots, I_n)$ of subsets $I_r \subseteq [n] := \{1,\ldots,n\}$ such that, for $i \in [n]$, if $i \in I_i$ then $I_{i+1}= (I_i \setminus \{i \}) \cup \{j\}$, for some $j \in [n]$; and if $i \in I_i$ then $I_{i+1} = I_i$. (Here the indices are taken modulo $n$.) In particular, we have $|I_1| = \cdots = |I_n|$.
\end{definition}

It is easy to see the bijection between necklaces and decorated permutations. To go from a Grassmann necklace $\I$ to a decorated permutation $\pi^{:}=(\pi,col)$:

\begin{itemize}
 \item if $I_{i+1} = (I_i \backslash \{i\}) \cup \{j\}$, $j \not = i$, then $\pi(i)=j$,
 \item if $I_{i+1} = I_i$ and $i \not \in I_i$ then $\pi(i)=i, col(i)=1$,
 \item if $I_{i+1} = I_i$ and $i \in I_i$ then $\pi(i)=i, col(i)=-1$.
\end{itemize}

To go from a decorated permutation $\pi^{:}=(\pi,col)$ to a Grassmann necklace $\I$,
$$I_r = \{ i \in [n] | i<_r \pi^{-1}(i) \textbf{ or } (\pi(i)=i \textbf{ and } col(i)=-1) \},$$
where we define the \newword{cyclically shifted ordering} $<_t$ on $[n]$ by the total order $t <_t t+1 <_t \cdots <_t n <_t 1 \cdots <_t t-1$.

Let us look at an example. For a decorated permutation $\pi^{:}$ with $\pi = 61482735$, we get $I_1 = \{1,2,3,5\},I_2 = \{2,3,5,6\},I_3 = \{3,5,6,1\},I_4 = \{4,5,6,1\},I_5 = \{5,6,8,1\}, I_6 = \{6,8,1,2\}, I_7 = \{7,8,1,2\}$ and $I_8 = \{8,1,2,3\}$.






A theorem from \cite{MR2834184} allows us to view positroids as objects coming from Grassmann necklaces. For $I,J \in {[n] \choose k}$, where
$$I=\{i_1, \cdots, i_k \}, i_1 <_t i_2 \cdots <_t i_k$$ and
$$J=\{j_1, \cdots, j_k \}, j_1 <_t j_2 \cdots <_t j_k,$$
we set 
$$I \leq_t J \text{ if and only if } i_1 \leq_t j_1, \cdots, i_k \leq_t j_k.$$

\begin{theorem}[\cite{MR2834184}]
\label{thm:mainO}
$\M$ is a positroid if and only if for some Grassmann necklace $(I_1, \cdots, I_n)$, the following condition holds: $H \in \M$ if and only if $H \geq_t I_t$ for any $t\in [n]$.
\end{theorem}

From the definition of positroids, it is easy to see that a restriction of a positroid is also a positroid. Theorem 19 of \cite{MR2834184} implies that the dual of a positroid is also a positroid. Therefore, the class of positroids is closed under contraction, restriction and dual (This was proven separately in \cite{war}). In this paper, we will describe contraction and restriction of a positroid in terms of its decorated permutation (Theorem 19 of \cite{MR2834184} implies that the decorated permutation of the dual positroid is obtained by taking the inverse permutation).

Due to the cyclic nature of positroids, we will think of the base set $[n]$ as elements of $\Z_n$. Draw a circle, put integers $1$ to $n$ in a clockwise order. $(a,b)$ means integers between $a$ and $b$ on the circle read in clockwise order. $(a,b], [a,b), [a,b]$ are defined similarly. For $a \in [n]$, $max_a(T)$ denotes the maximal element of the set $T$ under the ordering $<_a$.

We end the section with a tool that we will use later.
\begin{lemma}
\label{lem:interval}
Let $A,B$ be $k$-element subsets of $[n]$ such that for any $a \in A \setminus B$ and $b \in B \setminus A$, we have $a \leq_t b$. Then $A \leq_t B$.
\end{lemma}
\begin{proof}
We write $A$ as $\{a_1 <_t \cdots <_t a_k\}$ and $B$ as $\{b_1 <_t \cdots <_t b_k\}$. When $k=1$, the statement is obvious. For the sake of induction, assume the claim is true for $k-1$. By setting $A' = \{a_2,\ldots,a_k\}$ and $B' = \{b_2,\ldots,b_k\}$, we get $A' \leq_t B'$ by induction hypothesis. Since we also have $a_1 \leq_t b_1$, we end up with $A \leq_t B$.
\end{proof}



\medskip

\section{Contraction in terms of Grassmann necklaces}

In this section, we will describe the contraction of positroids in terms of Grassmann necklaces. Let us pick a positroid $\M \subseteq {[n] \choose k}$ and its corresponding Grassmann necklace $\I_{\M}=(I_1,\cdots,I_n)$. Denote the associated decorated permutation as $\pi^{:}=(\pi,col)$. When $j$ is a fixed point of $\pi$, $j$ is either a coloop (when $col(j)=-1$) or a loop (when $col(j)=+1)$ of $\M$. In the previous case, the decorated permutation of the new positroid can be obtained by switching the color of $j$, and in the latter case, we get an empty positroid with decorated permutation $\pi = 1\cdots n$, $col(i) = +1$ for all $i \in [n]$.

Therefore, we only have to study the case when $j$ is not a fixed point of $\pi$, which implies that $j \in I_i$ for $i \in (\pi^{-1}(j),j]$.

Let $\M'$ be the positroid $\{I \in \M | j \in I\}$ (if we switch the sign of $col(j)$ for the decorated permutation corresponding to $\M'$, we would get a decorated permutation of $\M / \{j\}$). Denote the Grassmann necklace of $\M'$ as $(K_1,\cdots,K_n)$.

\begin{lemma}
Choose any $a \in [n]$, such that $j \not \in I_a$. Then we have $(I_a \setminus \{max_a(I_a \setminus I_j )\}) \cup \{j\} \in \M'$.
\end{lemma}

\begin{proof}
We will write $(I_a \setminus \{max_a(I_a \setminus I_j )\}) \cup \{j\} \in \M'$ as $H$. We use $h_1 <_a \cdots <_a h_k$ to denote the elements of $I_a$. From the definition of Grassmann necklace, we have $I_j \setminus I_a \subseteq [j,a)$ and $I_a \setminus I_j \subseteq [a,j)$. Let $h_i$ be the element $max_a(I_a \setminus I_j)$, which implies that $h_{i+1},\ldots,h_k \in I_j$.

Since $H = I_a \setminus \{h_i\} \cup \{j\}$, we have $H \setminus I_j \subseteq [a,h_i)$ and $I_j \setminus H \subseteq (j,a)$. This combined with Lemma~\ref{lem:interval} gives us $I_j \leq_t H$ for $t \in [h_i,j]$. Again using $H = I_a \setminus \{h_i\} \cup \{j\}$, we get $I_a \leq_t H$ for $t \in (j,h_i]$. Combining the two results, we can conclude that $I_t \leq_t H$ for all $t \in [n]$, netting us with $H \in \M$ via Theorem~\ref{thm:mainO}. 





\end{proof}

\begin{proposition}
If $j \in I_a$, we have $K_a = I_a$. If not, then we have $K_a = (I_a \setminus \{max_a(I_a \setminus I_j )\}) \cup \{j\}$.
\end{proposition}

\begin{proof}
When $j \in I_a$, it is obvious. Hence let us look at the case when $j \not \in I_a$. Using the same notation as in the proof of the previous lemma, set $q$ such that $h_q <_a j <_a h_{q+1}$. Then we can write the elements of $I_a$ and $H$ in increasing order from left to right with respect to $<_a$ as:
$$I_a = \{h_1,\cdots,h_{i-1},h_i,\cdots,h_q,\hat{j},h_{q+1},\cdots,h_k\},$$
$$H = \{h_1,\cdots,h_{i-1},\hat{h_i},\cdots,h_q,j,h_{q+1},\cdots,h_k\},$$
where $\hat{}$ means the corresponding element is not contained in the set.

From $j \in K_a$ and the relationship $I_a \leq_a K_a \leq_a H$, we can see that
$$K_a = \{h_1,\cdots,h_{i-1},s_i,\cdots,s_{q-1},j,h_{q+1},\cdots,h_k\},$$
for some elements $s_i,\ldots,s_{q-1}$. Using $K_j = I_j$, $\{h_{i+1},\ldots,h_k\} \subseteq I_j$ and $K_j \setminus K_a \subseteq [j,a)$, we end up with $s_i = h_{i+1},\ldots,s_{q-1} = h_q$.


\end{proof}

Let us look at an example. When $\pi=[6,1,4,8,2,7,3,5]$ and $j=3$:
$$ I_1 = \{1,2,3,5\}, K_1 = \{1,2,3,5\},  $$
$$ I_2 = \{2,3,5,6\}, K_2 = \{2,3,5,6\}, $$
$$ I_3 = \{3,5,6,1\}, K_3 = \{3,5,6,1\}, $$
$$ I_4 = \{4,5,6,1\}, K_4 = \{5,6,1,3\}, $$
$$ I_5 = \{5,6,8,1\}, K_5 = \{5,6,1,3\}, $$
$$ I_6 = \{6,8,1,2\}, K_6 = \{6,8,1,3\}, $$
$$ I_7 = \{7,8,1,2\}, K_7 = \{7,8,1,3\}, $$
$$ I_8 = \{8,1,2,3\}, K_8 = \{8,1,2,3\}. $$

$K_1,K_2,K_3,K_8$ equals $I_1,I_2,I_3,I_8$ respectively, since they contain $3$. For $K_4$, $max_4(\{4\})=4$ gives us $K_4=(I_4 \setminus \{4\}\cup \{3\}) =\{3,5,6,1\}$. 

\section{Contraction in terms of Decorated permutations}
In this section, we will describe how the decorated permutation of $\M':= \M/\{j\}$ would look like. The goal is to describe the decorated permutation of $\M'$ without having to compute all $K_a$'s. When $\pi(j) = j$, the problem is trivial, so let us assume that $\pi(j) \not = j$. We denote $\mu^{:}=(\mu,col')$ to be the decorated permutation of $\M'$. Throughout the section, we use
$$I_a \stackrel{\pi(a)}{\longrightarrow} I_{a+1}$$
to express 
$$I_{a+1} = \begin{array}{cc} I_a &  \textit{If   }  a \not \in I_a,

\\(I_a \setminus \{a\}) \cup \{\pi(a)\} & \textit{Otherwise.} \end{array} $$ 

We will also use
$$\begin{CD}
I_{a} \\
@VV\phi(a)\ V \\
K_{a} 
\end{CD}$$
to denote $K_a = (I_a \setminus \{\phi(a)\}) \cup \{j\}$, where 
$$\phi(a) := \begin{array}{cc} j & \textit{If   } j \in I_a

\\max_a(I_a \setminus I_j ) & \textit{Otherwise.} \end{array} $$

Now let us look at the following diagram concerning $I_a,K_a,I_{a+1},K_{a+1}$.

$$\begin{CD}
I_{a} @>\pi(a)>> I_{a+1}\\
@VV\phi(a)V @VV\phi(a+1)V\\
K_{a} @>\mu(a)>> K_{a+1}
\end{CD}$$

We will refer to this diagram as a \newword{square} at $a$. One of the following two cases is possible:
\begin{itemize}
\item $\pi(a)=\phi(a+1)$ and $\phi(a)=\mu(a),$
\item $\pi(a)=\mu(a)$ and $\phi(a+1)=\pi(a).$
\end{itemize}

\medskip

Let us try to analyze those squares, in order to describe $\mu(a)$ using only $\pi^{:}$.

\begin{itemize}

\item Case 1 : When $a=j$, we get $\mu(j)=j$ since $j \in K_i$ for all $i \in [n]$. 

$$\begin{CD}
I_{j} @>\pi(j)>> I_{j+1}\\
@VVjV @VV\phi(j+1)=\pi(j)V\\
K_{j} @>\mu(j)=j>> K_{j+1}
\end{CD}$$

\item Case 2 : When $a \in (\pi^{-1}(j),j)$, we get $\mu(a)=\pi(a)$ from the fact that $j \in I_t$ for all $t \in (\pi^{-1}(j),j]$.

$$\begin{CD}
I_{a} @>\pi(a)>> I_{a+1}          \\
@VV\phi(a)=jV @VV\phi(a+1)=jV             \\
K_{a} @>\mu(a)=\pi(a)>> K_{a+1}  
\end{CD}
$$

\item Case 3 :When $a = \pi^{-1}(j)$, we get $\mu(a) = \phi(\pi^{-1}(j))=\pi(j)$ from the fact that $j \in I_{\pi^{-1}(j)+1}$.

$$\begin{CD}
I_{\pi^{-1}(j)} @>j>> I_{\pi^{-1}(j)+1}\\
@VV\phi(\pi^{-1}(j))V @VVjV\\
K_{\pi^{-1}(j)} @>\mu(a)=\phi(\pi^{-1}(j))>> K_{\pi^{-1}(j)+1}
\end{CD}$$

\item When $a \in (j,\pi^{-1}(j))$, we have three cases:
\begin{itemize}
\item Case 4a : If $\phi(a)=a$, from $\phi(a)=a \not \in K_a$, we end up with $\mu(a)=a$.
$$\begin{CD}
I_{a} @>\pi(a)>> I_{a+1}          \\
@VV\phi(a)=a V @VV\phi(a+1)=\pi(a)V             \\
K_{a} @>\mu(a)=a>> K_{a+1}  
\end{CD}
$$

\item Case 4b : When $\phi(a) \not =a$ and $\pi(a) >_{a+1} j$, we get $\phi(a)=\phi(a+1)$. This is because we can't have $\phi(a+1)=\pi(a)$ due to $\pi(a) >_{a+1} j$.
$$\begin{CD}
I_{a} @>\pi(a)>> I_{a+1}          \\
@VV\phi(a)V @VV\phi(a+1)=\phi(a)V             \\
K_{a} @>\mu(a)=\pi(a)>> K_{a+1}  
\end{CD}
$$

\item Case 4c : When $\phi(a) \not =a$ and $\pi(a) <_{a+1} j$, we have $\phi(a)= max_{a+1}(\pi(a),\phi(a))$ and $\mu(a)=min_{a+1}(\pi(a),\phi(a))$. 
\end{itemize}

$$\begin{CD}
I_{a} @>\pi(a)>> I_{a+1}          \\
@VV\phi(a)V @VV\phi(a+1) =max_{a+1}(\pi(a),\phi(a)) V             \\
K_{a} @>\mu(a)=min_{a+1}(\pi(a),\phi(a)) >> K_{a+1}  
\end{CD}
$$

\end{itemize}

\medskip

Using what we have obtained so far, we can express each $\phi(a)$ without computing the $K_a$'s. We will first take a look at an example, when $\pi=[6,1,4,8,2,7,3,5]$ and $j=3$.

$$\begin{CD}
1235 @>6>> 2356 @>1>> 3561 @>4>> 4561 @>8>> 5681 @>2>> 6812 @>7>> 7812 @>3>> 8123 @>5>>         \\
@VV3V @VV3V @VV3V @VV4V @VV8V @VV2V @VV2V @VV3V              \\
1235 @>6>> 2356 @>1>> 3561 @>3>> 5613 @>4>> 5613 @>8>> 6813 @>7>> 7813 @>2>> 8123 @>5>>
\end{CD}
$$

Start with the square at $3$. $\phi(4)$ should equal $\pi(3)$. As we move to the right, $\phi(a+1)$ should always equal $\pi(a)$ or $\phi(a)$. In other words, we build $\mu$ by starting from the square at $j$, then moving to the right. We use $a$ to denote which square we are currently looking at.

We start out by setting $\mu^{:}$ to be $\pi^{:}$, and as we move along the squares, we change $\mu(a)$'s accordingly. The first square, square $j$, corresponds to case 1 above. We set $\mu(j)=j$ and $q = \phi(j+1) = \pi(j)$. As we move to the right, when we are considering a square at $a$, $q$ will be the record for $\phi(a)$ that was determined in the previous square.

When we are dealing with square at $a$, if $\mu(a)=\pi(a)$, we have $\phi(a+1) = \phi(a)$ and hence $q$ doesn't change. We change $\mu(a)$ to $q$ and $q$ to $\pi(a)$ when we are in one of the following cases: 

\begin{enumerate}
\item Case $4a$,
\item Case $4c$ and we have $q <_{a+1} \pi(a)$.
\end{enumerate}

That is, when $q=a$ or $q<_{a+1} \pi(a) <_{a+1} j$. We repeat this procedure until we end up at square $\pi^{-1}(j)$, which corresponds to case 3 above. 

\begin{algorithm}
\caption{ When $\pi(j) \not = j$ and $\M_{\mu^{:}} = \M_{\pi^{:}} / \{j\}$.  Obtaining $\mu^{:}=(\mu,col')$ from $\pi^{:}=(\pi,col)$.}
\begin{algorithmic}
\label{al:contraction} 
\STATE $\mu \Leftarrow \pi$
\STATE $col' \Leftarrow col$
\STATE $\mu(j) \Leftarrow j$
\STATE $col'(j) = 1$
\STATE $a \Leftarrow j+1$
\STATE $q \Leftarrow \pi(j)$

\WHILE{$\pi(a) \not = j$}
\IF{$q =a$ or $q <_{a+1} \pi(a) <_{a+1} j$} 
 \STATE $\mu(a) \Leftarrow q$
 \STATE $q \Leftarrow \pi(a)$
 \IF{ $\mu(a)=a$} 
 \STATE $col'(a) = 1$
 \ENDIF
\ENDIF
\STATE $a \Leftarrow a+1$
\ENDWHILE

\STATE $\mu(a) \Leftarrow q$

\end{algorithmic}
\end{algorithm}

When we do $a \Leftarrow a+1$, we are doing it modulo $n$. 

\begin{theorem}
Let $\pi^{:}$ be a decorated permutation. We contract the positroid corresponding to $\pi^{:}=(\pi,col)$ by $\{j\}$ to get a positroid corresponding to $\mu^{:}=(\mu,col')$. If $\pi(j) \not =j$, then $\mu^{:}$ is obtained by Algorithm~\ref{al:contraction}. If $\pi(j)=j$ and $col(j) = -1$, then $\mu^{:}=(\pi,col')$ where $col'(j)=1$ and $col'(i)=col(i)$ for all $i \not = j$. If $\pi(j)=j$ and $col(j)=+1$, we get $\mu = 123 \cdots n$ and $col(i)=1$ for all $i \in [n]$.
\end{theorem}

Let us try out the algorithm for $\pi=[6,1,4,8,2,7,3,5]$ and $j=3$.
\begin{enumerate}
\item $\mu = [6,1,4,8,2,7,3,5]$.
\item $\mu(3)=3, col'(3) = 1, q=4$. We start with $a=4, \mu=[6,1,3,8,2,7,3,5]$.
\item $a=4$ : From $q=a=4$, we get $\mu(4)=4, q=\pi(4)=8, col'(4)=1$. 
\item $a=5$ : Since $q=8 <_6 \pi(5)=2<_6 j=3$, we get $\mu(5)=8,q=2$. 
\item $a=6$ : Since $q \not =a$, $\pi(6)=7<_7 q=2$, we get $\mu(6)=7$ (unchanged) and $q=2$. 
\item $a=7$ : Since $\pi(7)=j=3$, we get $\mu(7)=2$ and we are done. 
\item We end with $\mu=[6,1,3,4,8,7,2,5]$.
\end{enumerate}

\section{Restriction}

In this section, we will work with restrictions. The proofs will be omitted due to duality. Let $\M \subseteq {[n] \choose k}$ be a positroid associated to a Grassmann necklace $\I_{\M}=(I_1,\cdots,I_n)$. Denote the corresponding decorated permutation as $\pi^{:}=(\pi,col)$. Define $\M'$ as $\{H \in \M | j \not \in H\}.$


Denote the Grassmann necklace of $\M'$ as $(K_1,\cdots,K_n)$, and its decorated permutation as $\mu^{:}=(\mu,col')$. 

\begin{lemma}
If $j \not \in I_a$, we have $K_a = I_a$. If not, then we have $K_a = (I_a \setminus \{j\}) \cup \{min_a(I_j \setminus I_a)\}$. 
\end{lemma}

Let us look at an example, where $\pi=[6,1,4,8,2,7,3,5]$ and $j=5$.
$$ I_1 = \{1,2,3,5\}, K_1 = \{1,2,3,6\} $$
$$ I_2 = \{2,3,5,6\}, K_2 = \{2,3,6,8\} $$
$$ I_3 = \{3,5,6,1\}, K_3 = \{3,6,8,1\} $$
$$ I_4 = \{4,5,6,1\}, K_4 = \{4,6,8,1\} $$
$$ I_5 = \{5,6,8,1\}, K_5 = \{6,8,1,2\} $$
$$ I_6 = \{6,8,1,2\}, K_6 = \{6,8,1,2\} $$
$$ I_7 = \{7,8,1,2\}, K_7 = \{7,8,1,2\} $$
$$ I_8 = \{8,1,2,3\}, K_8 = \{8,1,2,3\} $$

We have $K_6 = I_6, K_7 = I_7, K_8 = I_8$ since $I_6,I_7,I_8$ doesn't contain $j$. For $K_1$, since $min_1(I_5 \setminus I_1)=6$, we get $K_1 = \{1,2,3,6\}$. 

We are going to use diagrams as before. Horizontal arrows denote the same thing. Now the vertical arrow $I_a \stackrel{h}{\longrightarrow} K_a$ denotes $K_a = I_a \setminus \{j\} \cup \{h\}$. For each $a \in [n]$, denote such $h$ by $\chi(a)$. 

The following is an example of a diagram when $\pi=[6,1,4,8,2,7,3,5]$ and $j=5$.

$$\begin{CD}
1235 @>6>> 2356 @>1>> 3561 @>4>> 4561 @>8>> 5681 @>2>> 6812 @>7>> 7812 @>3>> 8123 @>5>>         \\
@VV6V @VV8V @VV8V @VV8V @VV2V @VV5V @VV5V @VV5V              \\
1236 @>8>> 2368 @>1>> 3681 @>4>> 4681 @>2>> 6812 @>5>> 6812 @>7>> 7812 @>3>> 8123 @>6>>
\end{CD}
$$

The main difference with the case of contractions is, we will be moving left starting at square at $j$. As we move to the left, we attach squares so that $\chi(a)$ equals $\chi(a+1)$ or $\pi(a)$.  

\begin{algorithm}
\caption{When $\pi(j) \not = j$ and $\M_{\mu^{:}} = \M_{\pi^{:}}|_{[n] \setminus \{j\}}$. Obtaining $\mu^{:}=(\mu,col')$ from $\pi^{:}=(\pi,col)$.}
\begin{algorithmic}
\label{al:restriction} 
\STATE $\mu \Leftarrow \pi$
\STATE $col' \Leftarrow col$
\STATE $\mu(j) \Leftarrow j$
\STATE $col'(j) = 1$
\STATE $a \Leftarrow j-1$
\STATE $q \Leftarrow \pi(j)$

\WHILE{$\pi(a) \not = j$}
\IF{$q =a$ or $q >_a \pi(a) >_{a} j$}
 \STATE $\mu(a) \Leftarrow q$
 \STATE $q \Leftarrow \pi(a)$
 \IF{$\mu(a)=a$}
 \STATE $col'(a) = 1$
 \ENDIF
\ENDIF
\STATE $a \Leftarrow a-1$
\ENDWHILE

\STATE $\mu(a) \Leftarrow q$

\end{algorithmic}
\end{algorithm}

\begin{theorem}
Let $\pi^{:}$ be a decorated permutation. We restrict the positroid corresponding to $\pi^{:}=(\pi,col)$ to $[n] \setminus \{j\}$, getting a positroid corresponding to $\mu^{:}=(\mu,col')$. If $\pi(j) \not = j$, then $\mu^{:}$ is obtained by Algorithm~\ref{al:restriction}. If $\pi(j)=j$ and $col(j) = +1$, then $\mu^{:}=\pi^{:}$. If $\pi(j)=j$ and $col(j)=-1$, we get $\mu = 123 \cdots n$ and $col(i)=1$ for all $i \in [n]$.
\end{theorem}

Let us try the algorithm for $\pi=[6,1,4,8,2,7,3,5]$ and $j=5$.
\begin{enumerate}
\item $\mu = [6,1,4,8,2,7,3,5]$.
\item $\mu(5)=5, col'(5)=1, q=2$. We start with $a=5, \mu=[6,1,4,8,5,7,3,5]$.
\item $a=4$ : From $q=2>_4 \pi(4)=8 >_4 j=5$, we get $\mu(4)=2, q=8$. 
\item $a=3$ : Since $q \not =a$ and $\pi(3)=4 <_3 j=5$, we get $\mu(3)=4$ (unchanged) and $q=8$. 
\item $a=2$ : Since $q \not =a$ and $\pi(2)=1 >_2 q=8$, we get $\mu(2)=1$ (unchanged) and $q=8$. 
\item $a=1$ : Since $q=8 >_1 \pi(1)=6 >_1 j=5$, we get $\mu(1)=8,q=6$. 
\item $a=8$ : Since $\pi(8)=j$, we get $\mu(8)=6$ and we are done.
\item We end with $\mu=[8,1,4,2,5,7,3,6]$.
\end{enumerate}





\bibliographystyle{plain}    
\bibliography{rs}        


\end{document}